\documentclass[11pt]{amsart}
\usepackage{}
\usepackage{amssymb}
\theoremstyle{plain}
\newtheorem{Thm}{Theorem}

\newtheorem{Lem}[Thm]{Lemma}
\newtheorem{Def}[Thm]{Definition}

\begin{document}

%begin Topmatter
\title[Convergence of Ricci flow to plane]
{Convergence of Ricci flow on $R^2$ to plane}

\author{Li Ma}
\address{Distinguished Professor, Department of mathematics \\
Henan Normal university \\
Xinxiang, 453007 \\
China} \email{nuslma@gmail.com}

\thanks{The research is partially supported by the National Natural Science
Foundation of China 10631020 and SRFDP 20090002110019}

\begin{abstract}
In this paper, we give a sufficient condition such that the Ricci
flow in $R^2$ exists globally and the flow converges at $t=\infty$
to the flat metric on $R^2$.


\end{abstract}

 \maketitle

\section{introduction}
In this short note, we are interested in the long-term behavior on
$R^2$ of conformally flat solutions to the Ricci flow equation on
$R^2$. Recall here that the Ricci flow equation for the
one-parameter family of metric $g(t)$ on $R^2$ is
\begin{equation}\label{ricci}
\partial_t g=-Rg,  \ \ in \ R^2.
\end{equation}
For these metrics $g(t)$, we take their forms as
$g(x,t)=e^{u(x,t)}g_E$, where $g_E$ is the standard Euclidean metric
on $R^2$. Then the Ricci flow equation becomes
\begin{equation}\label{scalar}
\partial_t e^{u}=\Delta u, \ \ in \ \ R^2,
\end{equation}
where $\Delta $ is the standard Laplacian operator of the flat
metric $g_E$ in $R^2$. The long-term existence of solutions of
(\ref{ricci}) or (\ref{scalar}) has been studied in \cite{DD}, ,
where it is shown that

\begin{Thm}\label{da}
 The solutions to (\ref{ricci}) with initial metric $g(0)=e^{u_0}g_E$ exist for all
$t\geq 0$ if and only if
\begin{equation}\label{pino}
\int_{R^2} e^{u_0}dx=\infty.
\end{equation}
\end{Thm}

The global behavior of the Ricci flow has been studied in \cite{W}.
To state one of her result, we recall two concepts of the metric
$g=g(t)$. One is below.

\begin{Def}\label{ap}
The \emph{aperture} of the metric $g$ on $R^2$ is defined as
$$ A(g)=\frac{1}{2\pi}\lim_{r\to\infty} \frac{L(\partial B_r)}{r}.
$$
 Here $B_r$ denotes the geodesic ball (or disc) of radius $r$
and $L(\partial B_r)$ is the length of the boundary of $\partial
B_r$. \end{Def}

The other is the Cheeger-Gromov convergence of the Ricci flow.

\begin{Def}\label{conv}
The Ricci flow $g(t)$ is said to have modified subsequence
convergence, if there exists a 1-parameter family of diffeomorphisms
$\{\phi(t)\}_{t_j\geq 0}$ such that for any sequence $t_j\to\infty$,
there exists a subsequence (denoted again by $t_j$ ) such that the
sequence $\phi_(t_j)^*g(t_j)$ converges uniformly on every compact
set as $t_j\to \infty$.
\end{Def}

Then we have the following result of L.F. Wu \cite{W}.

\begin{Thm}\label{wu}
Let $g(t) = e^{u(t)}g_E$ be a solution to (1.1) such that $g(0)=
e^{u_0}g_E$ is a complete metric with bounded curvature and $\nabla
u_0|$ is uniformly bounded on $R^2$. Then the Ricci flow has
modified subsequence convergence as $t_j\to \infty$ with the
limiting metric $g_\infty$ being complete metric on $R^2$;
furthermore, the limiting metric is flat if $A(g(0)>0$.
\end{Thm}
We point out that the diffeomorphisms $\phi(t_j)$ used in Theorem
\ref{wu} are of the special form
$$
\phi(t)(a,b)=(e^{\frac{-u(x_0,t)}{2}}a,e^{\frac{-u(x_0,t)}{2}}b)=(x_1,x_2)=x,
$$
where $x_0=(0,0)$. The important fact for these diffefeomorphisms is
that$$ |\nabla_{g(t)} f(x, t)|=|\nabla_{\phi(t)^*g(t)} f((a, b),
t)|.
$$ for any smooth function $f$ and $x=\phi(t)(a,b)$.

In the interesting paper \cite{IJ}, which motivates our work here,
the authors have proved the following.

\begin{Thm}\label{isen}
Suppose $g_0= e^{u_0}g_E$ has bounded curvature and $u_0$ is a
bounded smooth function on $R^2$. Then the Ricci flow $\partial_tg =
-Rg$ exists for all $t\geq 0$ and has modified subsequence
convergence to the flat metric in the $C^k$ topology of metrics on
compact domains in $R^2$ for each $k\geq 2$.
\end{Thm}

There is another formulation in dimension two. Since every complete
Riemannian manifold of dimension two is a one dimension K\"ahler
manifold, we can use the K\"ahler-Ricci flow formulation of the
Ricci flow on $R^2$. We shall consider the Ricci flow (\ref{ricci})
as the K\"ahler-Ricci flow by setting
$$
g_{i\bar{j}}={g_0}_{i\bar{j}}+\partial_i\partial_{\bar{j}}\phi,
$$
where $\phi=\phi(t)$ is the K\"ahler potential of the metric $g(t)$
relative to the metric $g_0$. Note that $$
g(0)_{i\bar{j}}={g_0}_{i\bar{j}}+\partial_i\partial_{\bar{j}}\phi_0,
$$
In this situation, the Ricci flow can be written as
\begin{equation}\label{kahler}
\partial_t\phi=\log \frac{{g_0}_{1\bar{1}}+\phi_{1\bar{1}}}{{g_0}_{1\bar{1}}}-f_0, \ \ \phi(0)=\phi_0,
\end{equation}
where $f_0$ is the potential function of the metric $g_0$ in the
sense that $R(g_0)=\Delta_{g_0}{f_0}$ in $R^2$. Here
$\Delta_{g_0}=g_0^{1\bar{1}}\partial_1\partial_{\bar{1}}$ in $R^2$,
which is the normalized Laplacian in K\"ahler geometry. Such a
potential function has been introduced by R.Hamilton in \cite{H}. We
remark that the initial data for the evolution equation
(\ref{kahler}) is $\phi(0)$ which is non-trivial. For the equivalent
of these two flows, one may see \cite{Ch}.

 Our result is below.
\begin{Thm}\label{ma}
Suppose $g_0= e^{u_0}g_E$ has bounded curvature $R_0$ with
(\ref{pino}) and $f_0$ is a bounded smooth function on $R^2$ such
that $\Delta_{g_0}f_0=R_0$. Then the Ricci flow $\partial_tg = -Rg$
with the initial metric $g_0$ exists for all $t\geq 0$ and has
modified subsequence convergence to the flat metric in the $C^k$
topology of metrics on compact domains in $R^2$ for each $k\geq 2$.
\end{Thm}

We remark that because of the assumption about the potential
function $f_0$, the initial metric $g_0$ is far from the cigar
metric \cite{M}. Here is the idea of the proof. We shall show that
the limit $f_\infty$ of $f(t_j)$ is a constant function. Because of
Theorem \ref{wu}, we need only show that
$R(g_\infty)=\Delta_{g_\infty}f_\infty=0$. The proof of Theorem
\ref{ma} will be given in section \ref{three}.

\section{maximum principle and the equivalence of flows (\ref{kahler}) and (\ref{ricci}) in dimension two}\label{two}

Fist we recall the maximum principle for the Ricci flow with bounded
curvature.  Given the Ricci flow $g(t)$ on $R^2$ with bounded
curvature, we have the following well-known maximum principle.

\begin{Lem}\label{MP1}Fix any $T>0$. If $w(x,t)$ is a bounded smooth solution to the heat
equation $$\partial_tw = \Delta_{g(t)}w, \ \ R^2\times (0,T]$$ with
the bounded initial data $w(x,0)$, then $|w(x,t)|\leq
\sup_{R^2}|w(x, 0)|$ for all $t\in(0,T]$.
\end{Lem}

We now consider  the equivalence of the flows (\ref{kahler}) and
(\ref{ricci}) in dimension two. We use the argument in \cite{Ch}
(see Lemma 4.1 there). If $g(t)$ is the Ricci flow in (\ref{ricci}),
we define
$$
u(x,t)=\int_0^t\log\frac{g_{1\bar{1}}(x,s)}{g_{1\bar{1}}(x,0)}ds-tf(0)
$$
and
$$
S_{1\bar{1}}(x,t)=g_{1\bar{1}}(x,t)-g_{1\bar{1}}(x,0)-u_{1\bar{1}}(x,t).
$$
Then by direct computation we have
$$
\frac{dS_{1\bar{1}}(x,t)}{dt}=0, \ \ S_{1\bar{1}}(x,0)=0.
$$
Hence $S_{1\bar{1}}(x,t)=0$ for all $t>0$ and
$$g_{1\bar{1}}(x,t)=g_{1\bar{1}}(x,0)+u_{1\bar{1}}(x,t).$$

If $u=u(x,t)$ is a solution to (\ref{kahler}), then it is clear that
$$g_{1\bar{1}}(x,t)=g_{1\bar{1}}(x,0)+u_{1\bar{1}}(x,t)$$
satisfies (\ref{ricci}).

\section{proof of Theorem \ref{ma}}\label{three}
The idea of the proof of Theorem \ref{ma} is similar to the argument
in \cite{BK} and \cite{IJ}, see also \cite{M}.

Let
$$
f=-\partial_t\phi.
$$
Then, taking the time derivative of (\ref{kahler}), we have
\begin{equation}\label{potential}
\partial_tf=\Delta_g f, \ \ f(0)=-\partial_t\phi(0)=f_0.
\end{equation}
By Lemma \ref{MP1} we know that $f$ is uniformly bounded in $R^2$.
 The important fact for us
is that
\begin{equation}\label{scalar}
\Delta_g f=R.
\end{equation}
See \cite{M} for a proof of this. if $f_0$ has some decay at space
infinity, one can can the same decay by the argument of Dai-Ma
\cite{DM}.

 It is well-known that $R$ is uniformly bounded in any
finite interval and $|f_t|$ and $|\nabla f|^2$ are bounded for each
$t\geq 0$ (via the use of $f(x,t)=f(x,0)+int_0^t R(x,s)ds$).

Our next task is to obtain a better control on $|\nabla f|$ as
$t\to\infty$. To get this, we let
$$
F(x,t)=t|\nabla f|^2+f^2.
$$
Then we have
$$
\partial_tF\leq \Delta_g F, \ \ in \ R^2.
$$
Using the maximum principle (Lemma \ref{MP1}), we know that
$$
\sup_{x\in R^2}|\nabla f(x,t)|^2\leq \frac{C}{1+t}
$$
for some uniform constant $C>0$. Once we have this bound, we can
follow the argument in Lemmata 8,9,and 10 in \cite{IJ} to conclude
that the curvature bounds that there are uniform constants $C_k$,
for any $k\geq 1$, such that
\begin{equation}\label{bounds}
\sup_{R^2}|\nabla^kR(x,t)|^2\leq \frac{C_k}{(1+t)^{k+2}}, \ t>0.
\end{equation}

We are now ready to complete the proof of Theorem \ref{ma} \emph{
Proof of Theorem \ref{ma}}. We shall use the modified convergence
sequence $g(t_j)$ in Theorem \ref{wu}. We need only show that the
limiting metric has flat curvature and this will be obtained by
showing that the limiting function $f_\infty$ of $f(x,t_j)$ is
constant. Since $f(x,t)$ is uniformly bounded by a constant $K> 0$
on $R^2\times [0,\infty)$, for the fixed $x_0=(0,0)\in R^2$ and for
any sequence $t_j\to\infty$, there exists a subsequence, still
denoted by $t_j$, such that $c=\lim f(x_0, t_j)$ exists. By the
construction of the metrics $g(t_j)$, for any compact subset
$K\subset R^2$, the limiting metric $g_\infty$ is equivalent to any
$\phi(t_j)^*g(t_j)$ for every large $t$; that is, there is a uniform
constant $C=C(K)>0$ such that
$$d_t(x, x_0)\leq C d_{g_\infty}(x, x_0)$$ for every $x\in K$, where $d_t(x, x_0)$ is the distance between $x$ and
$x_0$ in $\phi(t)^*g(t)$ and $d_{g_\infty}(x, x_0)$ is the distance
of the limiting metric.. For $x\in K$, we can establish (for all
$t>1$),
$$ |f(x, t)-f(x_0, t)| \leq d_t(x, x_0) \sup_{x\in K}|\nabla
f(x, t)| \leq \frac{Cd_{g_\infty}(x, x_0)}{1+t},$$ where we have
used the fact that for $x=\phi(t)(a,b)$,$$ |\nabla_{g(t)} f(x,
t)|=|\nabla_{\phi(t)^*g(t)} f((a,b), t)|,
$$which is uniformly bounded. It follows that $f(x, t_j)$ is also convergent to $c$, which is
$f_\infty(x)=c$ as $t_j\to\infty$, and then
$\partial_1\partial_{\bar{1}} f(x,t_j)\to 0$. Then we have
$\Delta_{g_\infty}f_\infty=0=R_\infty$ and then
$\phi(t_j)^*g(t_j)\to g_\infty$ locally in $C^2$ with $g_\infty$ of
flat curvature. The $C^k$-convergence of $\phi(t_j)^*g(t_j)$ to this
flat limit then follows from the previous curvature estimates
obtained in (\ref{bounds}) (see Lemmata 8,9, and 10 in \cite{IJ}).
This completes the proof of Theorem \ref{ma}.

\end{document}